\documentclass[11pt]{amsart}
\usepackage{amscd,amssymb,amsmath}
\usepackage[arrow,matrix]{xy}
\usepackage{graphicx}
\usepackage{cite}
\usepackage{ifpdf}
\usepackage{tikz}

\newtheorem{thm}[subsection]{Theorem}

\newtheorem{defn}[subsection]{Definition}

\numberwithin{equation}{section} 

\newcommand{\M}{{\mathcal M}}
\newcommand{\A}{{\mathcal A}}

\newcommand{\bea}{\begin{eqnarray}}
\newcommand{\eea}{\end{eqnarray}}

\newcommand{\R}{\mathbb{R}}

\begin{document}
\title[Dynamics of two-dimensional evolution algebras]
{Dynamics of two-dimensional evolution algebras}

\author{U.A. Rozikov, Sh.N. Murodov}

 \address{U.\ A.\ Rozikov\\ Institute of mathematics,
Tashkent, Uzbekistan.} \email {rozikovu@yandex.ru}
 \address{Sh.\ N.\ Murodov\\ Institute of mathematics,
Tashkent, Uzbekistan.} \email {sherruriy$_-$1086@mail.ru}

\begin{abstract} Recently in \cite{clr} a notion of a chain of evolution algebras is introduced.
This chain is a dynamical system the state of which at each given time is an evolution algebra.
The sequence of matrices of the structural constants for this
chain of evolution algebras satisfies the Chapman-Kolmogorov equation.
In this paper we construct 25 distinct examples of chains of two-dimensional evolution algebras.
 For all of these 25 chains we study the behavior of the baric
property, the behavior of the set of absolute nilpotent elements and dynamics of the set of idempotent elements depending on the time.

{\it AMS classifications (2010):} 17D92; 37C99; 60J25.\\[2mm]

\end{abstract}
\maketitle

\section{Introduction} \label{sec:intro}

 The concept of a dynamical system has its origins in Newtonian mechanics.
 There, as in other natural sciences and engineering disciplines, the evolution rule of
  a dynamical systems
 is given implicitly by a relation that gives the state of the system only a short time
 into the future. The relation is either a differential equation, a difference equation or another time scale.
 To determine the state for all future times requires iterating the relation many times each
  advancing time by a small step.

For simple dynamical systems, knowing the trajectory is often sufficient, but most dynamical
 systems are too complicated to be understood in terms of individual trajectories.

Following Kolmogorov, an approach based on the time evolution of initial distributions
 appears to be more appropriate (see e.g.\cite{AFS}).

 In \cite{t} a notion of evolution algebra is introduced. This
evolution algebra is defined as follows.
 Let $(E,\cdot)$ be an algebra over a field $K$. If it admits a
basis $e_1,e_2,\dots$, such that $e_i\cdot e_j=0$, if $i\ne j$ and $
e_i\cdot e_i=\sum_{k}a_{ik}e_k$, for any $i$, then this algebra is
called an {\it evolution algebra}.

The concept of evolution algebras lies between algebras and dynamical
systems. Algebraically, evolution algebras are non-associative
Banach algebra; dynamically, they represent discrete dynamical
systems. Evolution algebras have many connections with other
mathematical fields including graph theory, group theory, stochastic
processes, mathematical physics, etc.

In the book \cite{t}, the foundation of evolution algebra theory and
applications in non-Mendelian genetics and Markov chains are
developed.

In \cite{rt} the algebraic structures of function spaces defined by
graphs and state spaces equipped with Gibbs measures by associating
evolution algebras are studied. Results of \cite{rt} also allow  a
natural introduction of thermodynamics in studying of several
systems of biology, physics and mathematics by theory of evolution
algebras.

There exist several classes of non-associative algebras (baric,
evolution, Bernstein, train, stochastic, etc.), whose investigation
has provided a number of significant contributions to theoretical
population genetics. Such classes have been defined different times
by several authors, and all algebras belonging to these classes are
generally called "genetic"  \cite{e1,e2,e3,LR,ly,m,t,w}.

Recently in \cite{clr} a notion of a chain of evolution algebras is introduced.
 This chain is a dynamical system the state of which at each given time is an evolution algebra.
 The chain is defined by the sequence of matrices of the structural constants (of the evolution
 algebras considered in \cite{t}) which satisfies the Chapman-Kolmogorov equation.
 In \cite{clr} several examples (time
homogenous, time non-homogenous, periodic, etc.) of such chains are given.

In this paper we continue investigation of chain of evolution algebras, we in more detail study chains generated by two-dimensional evolution algebras.

The paper is organized as follows. In Section 2 we give main
definitions and facts related to a chain of evolution algebras.
In Section 3 we construct 25 distinct chains of two-dimensional evolution algebras. In Section 4 we study
behavior of the property to be baric for each 25 chains constructed in the section 3.
We show that some of the chains are baric any time, some of them are never baric.
For other chains of 25 ones which have (baric property transition)
we define a baric property controller function
and under some conditions on this controller we prove that the
chain is not baric almost surely (with respect to Lebesgue
measure).  For each 25 chains of evolution algebras in Section 5 we study the behavior of the set of absolute nilpotent elements and in the last section we study dynamics of the set of idempotent elements depending on the time.

\section{Preliminaries}
Following \cite{clr} we consider a family
$$\left\{E^{[s,t]}:\ s,t \in \R,\ 0\leq s\leq t
\right\}$$ of $n$-dimensional evolution algebras over the field $\R$,
with basis $e_1,\dots,e_n$ and multiplication table
\begin{equation}\label{1}
 e_ie_i =\sum_{j=1}^na_{ij}^{[s,t]}e_j, \ \ i=1,\dots,n; \ \
e_ie_j =0,\ \ i\ne j.\end{equation} Here parameters $s,t$ are
considered as time.

Denote by
$\M^{[s,t]}=\left(a_{ij}^{[s,t]}\right)_{i,j=1,\dots,n}$-the matrix
of structural constants.

\begin{defn}\label{d1}\cite{clr} A family $\left\{E^{[s,t]}:\ s,t \in \R,\ 0\leq s\leq t
\right\}$ of $n$-dimensional evolution algebras over the field $\R$
is called a chain of evolution algebras (CEA) if the matrix
$\M^{[s,t]}$ of structural constants satisfies the
Chapman-Kolmogorov equation
\begin{equation}\label{2}
\M^{[s,t]}=\M^{[s,\tau]}\M^{[\tau,t]}, \ \ \mbox{for any} \ \
s<\tau<t.
\end{equation}
\end{defn}

\begin{defn}\label{d2} A CEA is called a time-homogenous CEA if
the matrix $\M^{[s,t]}$ depends only on $t-s$. In this case we write
$\M^{[t-s]}$.
\end{defn}

\begin{defn}\label{dp} A CEA is called periodic if its matrix $\M^{[s,t]}$
is periodic with respect to at least one of the variables  $s$, $t$,
i.e. (periodicity with respect to $t$) $\M^{[s,t+P]}=\M^{[s,t]}$ for
all values of $t$. The constant $P$ is called the period, and is
required to be nonzero.
\end{defn}

Let $\left\{\M^{[s,t]}, \ \ 0\leq s\leq t\right\}$ be a family of
stochastic matrices which satisfies the equation (\ref{2}), then it
defines a Markov process. Thus we have

\begin{thm}\label{t1}\cite{clr} For each Markov process, there is a CEA whose
structural constants are transition probabilities of the process,
and whose generator set (basis) is the state space of the Markov
process.
\end{thm}

A {\it character} for an algebra $A$ is a nonzero multiplicative
linear form on $A$, that is, a nonzero algebra homomorphism from $A$
to $\R$ \cite{ly}. Not every algebra admits a character. For
example, an algebra with the zero multiplication has no character.

\begin{defn}\label{d3} A pair $(A, \sigma)$ consisting of an algebra $A$ and a
character $\sigma$ on $A$ is called a {\it baric algebra}. The
homomorphism $\sigma$ is called the weight (or baric) function of
$A$ and $\sigma(x)$ the weight (baric value) of $x$.
\end{defn}
In \cite{ly} for the evolution algebra of a free population it is
proven that there is a character $\sigma(x)=\sum_i x_i$, therefore
that algebra is baric. But the evolution algebra $E$ introduced in
\cite{t} is not baric, in general. The following theorem gives a
criterion for an evolution algebra $E$ to be baric.

\begin{thm}\label{t2}\cite{clr} An $n$-dimensional evolution algebra $E$, over the field $\R$,
is baric if and only if there is a column
$\left(a_{1i_0},\dots,a_{ni_0}\right)^T$ of its structural constants
matrix $\M=\left(a_{ij}\right)_{i,j=1,\dots,n}$, such that
$a_{i_0i_0}\ne 0$ and $a_{ii_0}=0$, for all $i\ne i_0$. Moreover,
the corresponding weight function is $\sigma(x)=a_{i_0i_0}x_{i_0}$.
\end{thm}

In \cite{clr} several concrete examples of CEAs are given and their time-dynamics are studied.
In this paper we continue investigation of CEAs, in more detail we study CEAs generated by two-dimensional evolution algebras.

\section{Construction of chains of two-dimensional EAs}

In this section we shall construct several chains of two-dimensional evolution algebras, which were not considered in \cite{clr}.

To construct a chain of two-dimensional evolution algebra one has to solve equation (\ref{2}) for $2\times 2$ matrix
$\M^{[s,t]}$. This equation gives the following system of functional equations (with four unknown functions):
\begin{equation}\label{fe}\begin{array}{llll}
a_{11}^{[s,t]}=a_{11}^{[s,\tau]}a_{11}^{[\tau,t]}+a_{12}^{[s,\tau]}a_{21}^{[\tau,t]},\\[2mm]
a_{12}^{[s,t]}=a_{11}^{[s,\tau]}a_{12}^{[\tau,t]}+a_{12}^{[s,\tau]}a_{22}^{[\tau,t]},\\[2mm]
a_{21}^{[s,t]}=a_{21}^{[s,\tau]}a_{11}^{[\tau,t]}+a_{22}^{[s,\tau]}a_{21}^{[\tau,t]},\\[2mm]
a_{22}^{[s,t]}=a_{21}^{[s,\tau]}a_{12}^{[\tau,t]}+a_{22}^{[s,\tau]}a_{22}^{[\tau,t]}.\\
 \end{array}
 \end{equation}

 But the analysis of the system (\ref{fe}) is difficult. We shall consider several cases where the system is solvable:

{\bf Case} 1. $a_{11}^{[s,t]}=a_{22}^{[s,t]}=\alpha(s,t)$, $a_{12}^{[s,t]}=a_{21}^{[s,t]}=\beta(s,t)$.  Then equation (\ref{fe}) is reduced to
$$\alpha(s,t)=\alpha(s,\tau)\alpha(\tau,t)+\beta(s,\tau)\beta(\tau,t);$$
$$\beta(s,t)=\alpha(s,\tau)\beta(\tau,t)+\beta(s,\tau)\alpha(\tau,t).$$
Denote
\begin{equation}\label{fv}
f(s,t)=\alpha(s,t)+\beta(s,t), \ \ \varphi(s,t)=\alpha(s,t)-\beta(s,t),
\end{equation} then
the last system of functional equations can be written as
$$f(s,t)=f(s,\tau)f(\tau,t), \ \
\varphi(s,t)=\varphi(s,\tau)\varphi(\tau,t), \ s\leq \tau\leq t.$$ Both these equations are
known as Cantor's second equation which has
very rich family of solutions:

a) $f(s,t)\equiv 0$;

b) $f(s,t)={\Phi(t)\over \Phi(s)}$, where $\Phi$ is an arbitrary
function with $\Phi(s)\ne 0$;

c) $$f(s,t)=\left\{\begin{array}{ll}
1, \ \ \mbox{if}\ \ s\leq t<a,\\[2mm]
0, \ \ \mbox{if} \ \ t\geq a.\\
\end{array}\right. \ \ \mbox{where} \ \ a>0.$$
Similarly,

a') $\varphi(s,t)\equiv 0$;

b') $\varphi(s,t)={\Psi(t)\over \Psi(s)}$, where $\Psi$ is an arbitrary
function with $\Psi(s)\ne 0$;

c') $$\varphi(s,t)=\left\{\begin{array}{ll}
1, \ \ \mbox{if}\ \ s\leq t<b,\\[2mm]
0, \ \ \mbox{if} \ \ t\geq b.\\
\end{array}\right. \ \ \mbox{where} \ \ b>0.$$
Substituting these solutions into (\ref{fv}) and finding $\alpha(s,t)$ and $\beta(s,t)$ we get the following matrices

 $$\M_0^{[s,t]}=\left(\begin{array}{cc}
0 & 0\\
0 & 0
\end{array}\right);\ \ \M_1^{[s,t]}={1\over 2}\left(\begin{array}{cc}
{\Psi(t)\over \Psi(s)} & -{\Psi(t)\over \Psi(s)}\\[2mm]
-{\Psi(t)\over \Psi(s)} & {\Psi(t)\over \Psi(s)}
\end{array}\right);$$

$$\M_2^{[s,t]}={1\over 2}\left\{\begin{array}{ll}
\left(\begin{array}{cc}
1 & -1\\
-1 & 1
\end{array}\right), \ \ \mbox{if} \ \ s\leq t<b;\\[4mm]
\left(\begin{array}{cc}
0 & 0\\[2mm]
0 & 0
\end{array}\right),\ \ \mbox{if} \ \ t\geq b\\
\end{array}\right.;$$
$$\M_3^{[s,t]}={1\over 2}\left(\begin{array}{cc}
{\Phi(t)\over \Phi(s)} & {\Phi(t)\over \Phi(s)}\\[2mm]
{\Phi(t)\over \Phi(s)} & {\Phi(t)\over \Phi(s)}
\end{array}\right).$$

$$\M_4^{[s,t]}={1\over 2}\left(\begin{array}{cc}
{\Phi(t)\over \Phi(s)}+{\Psi(t)\over \Psi(s)} & {\Phi(t)\over \Phi(s)}-{\Psi(t)\over \Psi(s)}\\[2mm]
{\Phi(t)\over \Phi(s)}-{\Psi(t)\over \Psi(s)} & {\Phi(t)\over \Phi(s)}+{\Psi(t)\over \Psi(s)}
\end{array}\right).$$

$$\M_5^{[s,t]}={1\over 2}\left\{\begin{array}{ll}
\left(\begin{array}{cc}
{\Phi(t)\over \Phi(s)}+1 & {\Phi(t)\over \Phi(s)}-1\\
{\Phi(t)\over \Phi(s)}-1 & {\Phi(t)\over \Phi(s)}+1
\end{array}\right), \ \ \mbox{if} \ \ s\leq t<b;\\[4mm]
\left(\begin{array}{cc}
{\Phi(t)\over \Phi(s)} & {\Phi(t)\over \Phi(s)}\\[2mm]
{\Phi(t)\over \Phi(s)} & {\Phi(t)\over \Phi(s)}
\end{array}\right),\ \ \mbox{if} \ \ t\geq b\\
\end{array}\right.;$$

$$\M_6^{[s,t]}={1\over 2}\left\{\begin{array}{ll}
\left(\begin{array}{cc}
1 & 1\\
1 & 1
\end{array}\right), \ \ \mbox{if} \ \ s\leq t<a;\\[4mm]
\left(\begin{array}{cc}
0 & 0\\[2mm]
0 & 0
\end{array}\right),\ \ \mbox{if} \ \ t\geq a\\
\end{array}\right.;$$

$$\M_7^{[s,t]}={1\over 2}\left\{\begin{array}{ll}
\left(\begin{array}{cc}
1+{\Psi(t)\over \Psi(s)} & 1-{\Psi(t)\over \Psi(s)}\\
1-{\Psi(t)\over \Psi(s)} & 1+{\Psi(t)\over \Psi(s)}
\end{array}\right), \ \ \mbox{if} \ \ s\leq t<a;\\[4mm]
\left(\begin{array}{cc}
{\Psi(t)\over \Psi(s)} & -{\Psi(t)\over \Psi(s)}\\[2mm]
-{\Psi(t)\over \Psi(s)} & {\Psi(t)\over \Psi(s)}
\end{array}\right),\ \ \mbox{if} \ \ t\geq a\\
\end{array}\right.;$$

$$\M_8^{[s,t]}={1\over 2}\left\{\begin{array}{llll}
\left(\begin{array}{cc}
2 & 0\\
0 & 2
\end{array}\right), \ \ \mbox{if} \ \ s\leq t<\min\{a,b\};\\[4mm]
\left(\begin{array}{cc}
1 & -1\\
-1 & 1
\end{array}\right), \ \ \mbox{if} \ \  a\leq t<b, \, a<b;\\[4mm]
\left(\begin{array}{cc}
1 & 1\\
1 & 1
\end{array}\right), \ \ \mbox{if} \ \  b\leq t<a, \, a>b;\\[4mm]
\left(\begin{array}{cc}
0 & 0\\[2mm]
0 & 0
\end{array}\right),\ \ \mbox{if} \ \ t\geq \max\{a,b\}\\
\end{array}\right..$$

Thus in this case we have nine CEAs: $E_i^{[s,t]}$, $0\leq s\leq t$ which correspond
to $\M_i^{[s,t]}$, $i=0,1,\dots, 8$ listed above.

{\bf Case} 2. $a_{11}^{[s,t]}=a_{22}^{[s,t]}=\alpha(s,t)$, $a_{12}^{[s,t]}=-a_{21}^{[s,t]}=\beta(s,t)$.  Then equation (\ref{fe}) is reduced to
\begin{equation}\label{pp}\begin{array}{ll}
\alpha(s,t)=\alpha(s,\tau)\alpha(\tau,t)-\beta(s,\tau)\beta(\tau,t);\\[2mm]
\beta(s,t)=\alpha(s,\tau)\beta(\tau,t)+\beta(s,\tau)\alpha(\tau,t).\\
\end{array}
\end{equation}
It is easy to check that this system of functional equations has a solution
$$\alpha(s,t)=\cos(t-s), \ \ \beta(s,t)=\sin(t-s).$$
This gives the following matrix
\begin{equation}\label{m}
\M_9^{[s,t]}=\left(\begin{array}{cc}
\cos (t-s) &  \sin (t-s)\\
-\sin (t-s) &  \cos (t-s)
\end{array}\right).
\end{equation}

Note that this matrix defines a periodic CEA. But we do not know another solution
(with $\alpha\beta\ne 0$) of the system (\ref{pp}).

We denote by $E_9^{[s,t]}=E_9^{[t-s]}$ the CEA which corresponds to $\M_9^{[s,t]}$.

{\bf Case} 3. $a_{11}^{[s,t]}=a_{21}^{[s,t]}=\alpha(s,t)$,
$a_{12}^{[s,t]}=a_{22}^{[s,t]}=\beta(s,t)$. In this case the equation (\ref{fe}) is reduced to
$$\alpha(s,t)=\alpha(\tau,t)(\alpha(s,\tau)+\beta(s,\tau));$$
$$\beta(s,t)=\beta(\tau,t)(\alpha(s,\tau)+\beta(s,\tau)).$$
Denote
\begin{equation}\label{gd}
\gamma(s,t)=\alpha(s,t)+\beta(s,t), \ \ \delta(s,t)=\alpha(s,t)-\beta(s,t),
\end{equation} then
the last system of functional equations can be written as
$$\gamma(s,t)=\gamma(s,\tau)\gamma(\tau,t), \ \
\delta(s,t)=\gamma(s,\tau)\delta(\tau,t), \ s\leq \tau\leq t.$$ The first equation is
Cantor's second equation which has
solutions:

a) $\gamma(s,t)\equiv 0$;

b) $\gamma(s,t)={h(t)\over h(s)}$, where $h$ is an arbitrary
function with $h(s)\ne 0$;

c) $$\gamma(s,t)=\left\{\begin{array}{ll}
1, \ \ \mbox{if}\ \ s\leq t<a,\\[2mm]
0, \ \ \mbox{if} \ \ t\geq a.\\
\end{array}\right. \ \ \mbox{where} \ \ a>0.$$
Using these solution from the second equation we find $\delta$:

a') $\delta(s,t)\equiv 0$;

b') $\delta(s,t)={g(t)\over h(s)}$, where $g$ is an arbitrary
function;

c') $$\delta(s,t)=\left\{\begin{array}{ll}
\psi(t), \ \ \mbox{if}\ \ s\leq t<a,\\[2mm]
0, \ \ \mbox{if} \ \ t\geq a.\\
\end{array}\right. \ \ \mbox{where} \ \ a>0, \, \psi(t)\, \mbox{is an arbitrary function}.$$

Substituting these solutions into (\ref{gd}) we get the following (non-zero) matrices

$$\M_{10}^{[s,t]}={1\over 2}\left(\begin{array}{cc}
{h(t)+g(t)\over h(s)} & {h(t)-g(t)\over h(s)}\\[2mm]
{h(t)+g(t)\over h(s)} & {h(t)-g(t)\over h(s)}\\[2mm]
\end{array}\right).$$

$$\M_{11}^{[s,t]}={1\over 2}\left\{\begin{array}{ll}
\left(\begin{array}{cc}
1+\psi(t) & 1-\psi(t)\\
1+\psi(t) & 1-\psi(t)\\
\end{array}\right), \ \ \mbox{if} \ \ s\leq t<a;\\[4mm]
\left(\begin{array}{cc}
0 & 0\\[2mm]
0 & 0
\end{array}\right),\ \ \mbox{if} \ \ t\geq a\\
\end{array}\right..$$

Thus in this case we have two new CEAs: $E_i^{[s,t]}$, $0\leq s\leq t$ which correspond
to $\M_i^{[s,t]}$, $i=10,11$ listed above.

{\bf Case} 4. $a_{11}^{[s,t]}=a_{12}^{[s,t]}=\alpha(s,t)$,
$a_{21}^{[s,t]}=a_{22}^{[s,t]}=\beta(s,t)$.
In this case the equation (\ref{fe}) is reduced to
$$\alpha(s,t)=\alpha(s,\tau)(\alpha(\tau,t)+\beta(\tau,t));$$
$$\beta(s,t)=\beta(s,\tau)(\alpha(\tau,t)+\beta(\tau,t)).$$
Denote
\begin{equation}\label{gd1}
\gamma(s,t)=\alpha(s,t)+\beta(s,t), \ \ \delta(s,t)=\alpha(s,t)-\beta(s,t),
\end{equation} then
the last system of functional equations can be written as
$$\gamma(s,t)=\gamma(s,\tau)\gamma(\tau,t), \ \
\delta(s,t)=\delta(s,\tau)\gamma(\tau,t), \ s\leq \tau\leq t.$$
The first equation has
solutions:

a) $\gamma(s,t)\equiv 0$;

b) $\gamma(s,t)={h(t)\over h(s)}$, where $h$ is an arbitrary
function with $h(s)\ne 0$;

c) $$\gamma(s,t)=\left\{\begin{array}{ll}
1, \ \ \mbox{if}\ \ s\leq t<a,\\[2mm]
0, \ \ \mbox{if} \ \ t\geq a.\\
\end{array}\right. \ \ \mbox{where} \ \ a>0.$$
Using these solution from the second equation we find $\delta$:

a') $\delta(s,t)\equiv 0$;

b') $\delta(s,t)={g(s)h(t)}$, where $g$ is an arbitrary
function;

c') $$\delta(s,t)=\left\{\begin{array}{ll}
\psi(s), \ \ \mbox{if}\ \ s\leq t<a,\\[2mm]
0, \ \ \mbox{if} \ \ t\geq a.\\
\end{array}\right. \ \ \mbox{where} \ \ a>0, \, \psi(t)\, \mbox{is an arbitrary function}.$$
Consequently, we get the following (non-zero) new matrices

$$\M_{12}^{[s,t]}={1\over 2}\left(\begin{array}{cc}
h(t)\left({1\over h(s)}+g(s)\right) & h(t)\left({1\over h(s)}+g(s)\right)\\[2mm]
h(t)\left({1\over h(s)}-g(s)\right) & h(t)\left({1\over h(s)}-g(s)\right)\\[2mm]
\end{array}\right);$$

$$\M_{13}^{[s,t]}={1\over 2}\left\{\begin{array}{ll}
\left(\begin{array}{cc}
1+\psi(s) & 1+\psi(s)\\
1-\psi(s) & 1-\psi(s)\\
\end{array}\right), \ \ \mbox{if} \ \ s\leq t<a;\\[4mm]
\left(\begin{array}{cc}
0 & 0\\[2mm]
0 & 0
\end{array}\right),\ \ \mbox{if} \ \ t\geq a\\
\end{array}\right..$$

Hence in this case we have two new CEAs: $E_i^{[s,t]}$, $0\leq s\leq t$ which correspond
to $\M_i^{[s,t]}$, $i=12,13$ listed above.

{\bf Case} 5. $a_{12}^{[s,t]}\equiv 0$, (the case $a_{21}^{[s,t]}\equiv 0$ is similar) then the system (\ref{fe}) reduced to following

\begin{equation}\label{fe0}\begin{array}{lll}
a_{11}^{[s,t]}=a_{11}^{[s,\tau]}a_{11}^{[\tau,t]},\\[2mm]
a_{21}^{[s,t]}=a_{21}^{[s,\tau]}a_{11}^{[\tau,t]}+a_{22}^{[s,\tau]}a_{21}^{[\tau,t]},\\[2mm]
a_{22}^{[s,t]}=a_{22}^{[s,\tau]}a_{22}^{[\tau,t]}\\
 \end{array}
 \end{equation}
 The first and the third equations of the system (\ref{fe0}) are Cantor's second equations.
 Substituting solutions of these equations into the second equation of the system (\ref{fe0}) we find
 function $a_{21}^{[s,t]}$. Note that in many cases the second equation of the system (\ref{fe0}) will be reduced
 to Cantor's  first equation:
 $$\gamma(s,t)=\gamma(s,\tau)+\gamma(\tau,t),$$
 which also has
very rich family of solutions: $\gamma(s,t)=\Psi(t)-\Psi(s)$, where
$\Psi$ is an arbitrary function.
 Thus solving the system (\ref{fe0}) we obtain the following new matrices:

 $$ \M_{14}^{[s,t]}=\left(\begin{array}{cc}
{\Phi(t)\over \Phi(s)} & 0\\[2mm]
\Phi(t)\psi(s) &0\\[2mm]
\end{array}\right);$$

$$\M_{15}^{[s,t]}=\left\{\begin{array}{ll}
\left(\begin{array}{cc}
1 & 0\\
\psi(s) & 0
\end{array}\right), \ \ \mbox{if} \ \ s\leq t<a;\\[4mm]
\left(\begin{array}{cc}
0 & 0\\[2mm]
0 & 0
\end{array}\right),\ \ \mbox{if} \ \ t\geq a\\
\end{array}\right.;$$

$$\M_{16}^{[s,t]}=\left(\begin{array}{cc}
0 & 0\\[2mm]
{g(t)\over \psi(s)} & {\psi(t)\over \psi(s)}\\[2mm]
\end{array}\right).$$

$$\M_{17}^{[s,t]}=\left(\begin{array}{cc}
{\Phi(t)\over \Phi(s)}  & 0\\[2mm]
{\Phi(t)\over \psi(s)}(g(t)-g(s)) & {\psi(t)\over \psi(s)}\\[2mm]
\end{array}\right).$$

$$\M_{18}^{[s,t]}=\left\{\begin{array}{ll}
\left(\begin{array}{cc}
1 & 0\\
{h(t)-h(s)\over\psi(s)} & {\psi(t)\over\psi(s)}\\
\end{array}\right), \ \ \mbox{if} \ \ s\leq t<a;\\[4mm]
\left(\begin{array}{cc}
0& 0\\
{h(t)\over \psi(s)} & {\psi(t)\over \psi(s)}\\
\end{array}\right),\ \ \mbox{if} \ \ t\geq a\\
\end{array}\right.;$$

$$\M_{19}^{[s,t]}=\left\{\begin{array}{ll}
\left(\begin{array}{cc}
0 & 0\\
h(t) & 1\\
\end{array}\right), \ \ \mbox{if} \ \ s\leq t<b;\\[4mm]
\left(\begin{array}{cc}
0& 0\\
0 & 0\\
\end{array}\right),\ \ \mbox{if} \ \ t\geq b\\
\end{array}\right.;$$

$$\M_{20}^{[s,t]}=\left\{\begin{array}{ll}
\left(\begin{array}{cc}
{\Phi(t)\over \Phi(s)} & 0\\
\Phi(t)(v(t)-v(s)) & 1\\
\end{array}\right), \ \ \mbox{if} \ \ s\leq t<b;\\[4mm]
\left(\begin{array}{cc}
{\Phi(t)\over \Phi(s)} & 0\\[2mm]
\Phi(t)w(s) & 0
\end{array}\right),\ \ \mbox{if} \ \ t\geq b\\
\end{array}\right.;$$

$$\M_{21}^{[s,t]}=\left\{\begin{array}{llll}
\left(\begin{array}{cc}
1 & 0\\
v(t)-v(s) & 1
\end{array}\right), \ \ \mbox{if} \ \ s\leq t<\min\{a,b\};\\[4mm]
\left(\begin{array}{cc}
1 & 0\\
v(s) & 0
\end{array}\right), \ \ \mbox{if} \ \  b\leq t<a, \, a>b;\\[4mm]
\left(\begin{array}{cc}
0 & 0\\
v(t) & 1
\end{array}\right), \ \ \mbox{if} \ \  a\leq t<b, \, a<b;\\[4mm]
\left(\begin{array}{cc}
0 & 0\\[2mm]
0 & 0
\end{array}\right),\ \ \mbox{if} \ \ t\geq \max\{a,b\}\\
\end{array}\right.;$$

 Hence, in this case we have eight  new CEAs: $E_i^{[s,t]}$, $0\leq s\leq t$ which correspond
to $\M_i^{[s,t]}$, $i=14,\dots, 21$ listed above.

{\bf Case} 6. $a_{22}^{[s,t]}\equiv 0$, (the case $a_{11}^{[s,t]}\equiv 0$ is similar) then the system (\ref{fe}) reduced to following
 \begin{equation}\label{fe2}\begin{array}{llll}
a_{11}^{[s,t]}=a_{11}^{[s,\tau]}a_{11}^{[\tau,t]}+a_{12}^{[s,\tau]}a_{21}^{[\tau,t]},\\[2mm]
a_{12}^{[s,t]}=a_{11}^{[s,\tau]}a_{12}^{[\tau,t]},\\[2mm]
a_{21}^{[s,t]}=a_{21}^{[s,\tau]}a_{11}^{[\tau,t]},\\[2mm]
0=a_{21}^{[s,\tau]}a_{12}^{[\tau,t]}.\\
 \end{array}
 \end{equation}
 Multiplying the second and the third equations and using the fourth one we get
 $a_{12}^{[s,t]}a_{21}^{[s,t]}=0$. Thus this case reduced to the Case 5, hence does not
 give any new CEA.

   {\bf Case} 7. Assume $a_{11}^{[s,t]}+a_{12}^{[s,t]}=a_{21}^{[s,t]}+a_{22}^{[s,t]}=1$. Denote
 $\alpha(s,t)=a_{11}^{[s,t]}, \ \ \beta(s,t)=a_{21}^{[s,t]}$ then from (\ref{fe}) we get
 $$\alpha(s,t)=\alpha(s,\tau)\alpha(\tau,t)+(1-\alpha(s,t))\beta(\tau,t);$$
$$\beta(s,t)=\beta(s,\tau)\alpha(\tau,t)+(1-\beta(s,\tau))\beta(\tau,t).$$
Putting
\begin{equation}\label{fv3}
\gamma(s,t)=\alpha(s,t)+\beta(s,t), \ \ \delta(s,t)=\alpha(s,t)-\beta(s,t),
\end{equation}
 we obtain
 \begin{equation}\label{a11}\begin{array}{ll}
  \gamma(s,t)=\gamma(s,\tau)\delta(\tau,t)+\gamma(\tau,t)-\delta(\tau,t);\\[2mm]
\delta(s,t)=\delta(s,\tau)\delta(\tau,t).
\end{array}\end{equation}
The second equation of (\ref{a11}) has the following solutions:

1) $\delta(s,t)\equiv 0$;

2) $\delta(s,t)={\theta(t)\over \theta(s)}$, with $\theta(t)\ne 0$;

3)  $$\delta(s,t)=\left\{\begin{array}{ll}
1, \ \ \mbox{if}\ \ s\leq t<a,\\[2mm]
0, \ \ \mbox{if} \ \ t\geq a.\\
\end{array}\right. \ \ \mbox{where} \ \ a>0.$$

Substituting these solutions into the first equation of (\ref{a11}) we get the following

1') $\gamma(s,t)=f(t)$, where $f$ is an arbitrary function;

2') For $\delta(s,t)={\theta(t)\over \theta(s)}$ we get
\begin{equation}\label{g}
\tilde{\gamma}(s,t)=\tilde{\gamma}(s,\tau)+\tilde{\gamma}(\tau,t)-{1\over\theta(\tau)},
\end{equation}
where $\tilde{\gamma}(s,t)={\gamma(s,t)\over \theta(t)}$.
We shall find a solution of the equation (\ref{g}) which has the form
$$\tilde{\gamma}(s,t)=\lambda\cdot u(t)-\mu\cdot u(s), \ \ \lambda, \mu\in \R,\, \lambda\ne \mu.$$
We do not know another kind of solutions of the equation (\ref{g}).

Substituting this function into equation (\ref{g}) we get $u(t)={1\over (\lambda-\mu)\theta(t)}.$
Consequently, we obtain
$$\gamma(s,t)={\lambda\over \lambda-\mu}-{\mu\theta(t)\over (\lambda-\mu)\theta(s)}.$$

3') In case of 3) we get the following equation
$$\gamma(s,t)=\left\{\begin{array}{ll}
\gamma(s,\tau)+\gamma(\tau,t)-1, \ \ \mbox{if} \ \ s\leq t<a,\\
\gamma(\tau,t), \ \ \mbox{if} \ \  t\geq a.
\end{array}\right.
$$
This equation has the following solution
$$\gamma(s,t)=\left\{\begin{array}{ll}
1, \ \ \mbox{if} \ \ s\leq t<a,\\
g(t), \ \ \mbox{if} \ \  t\geq a,
\end{array}\right.
$$
where $g(t)$ is an arbitrary function.

Using the obtained solutions $\gamma$ and $\delta$ by (\ref{fv3}) we get the following matrices

 $$ \M_{22}^{[s,t]}=\left(\begin{array}{cc}
f(t)& 1-f(t)\\[2mm]
f(t) &1-f(t)\\[2mm]
\end{array}\right);$$

$$ \M_{23}^{[s,t]}(\lambda,\mu)=\left(\begin{array}{cc}
1-{\lambda-2\mu\over 2(\lambda-\mu)}\left(1-{\theta(t)\over\theta(s)}\right)& {\lambda-2\mu\over 2(\lambda-\mu)}\left(1-{\theta(t)\over\theta(s)}\right)\\[2mm]
{\lambda\over 2(\lambda-\mu)}\left(1-{\theta(t)\over\theta(s)}\right)& 1-{\lambda\over 2(\lambda-\mu)}\left(1-{\theta(t)\over\theta(s)}\right) \\[2mm]
\end{array}\right);$$

$$\M_{24}^{[s,t]}=\left\{\begin{array}{ll}
\left(\begin{array}{cc}
1 & 0\\
0 & 1
\end{array}\right), \ \ \mbox{if} \ \ s\leq t<a;\\[4mm]
\left(\begin{array}{cc}
g(t) & 1-g(t)\\[2mm]
g(t) & 1-g(t)
\end{array}\right),\ \ \mbox{if} \ \ t\geq a\\
\end{array}\right..$$

 In this case we have three  new CEAs: $E_i^{[s,t]}$, $0\leq s\leq t$ which correspond
to $\M_i^{[s,t]}$, $i=22,23,24$ listed above.

\section{Baric property transition}

 In \cite{clr} a notion of property
transition for CEAs is defined. We recall the definitions:

\begin{defn}\label{d4} Assume a CEA, $E^{[s,t]}$, has a property, say $P$,
at pair of times $(s_0,t_0)$; one says that the CEA has $P$ property
transition if there is a pair $(s,t)\ne (s_0,t_0)$ at which the CEA
has no the property $P$.
\end{defn}

Denote
$$\mathcal T=\{(s,t): 0\leq s\leq t\};$$
$$\mathcal T_P=\{(s,t)\in \mathcal T: E^{[s,t]} \ \ \mbox{has property} \  P \};$$
$$\mathcal T_P^0=\mathcal T\setminus \mathcal T_P=\{(s,t)\in \mathcal T: E^{[s,t]} \ \ \mbox{has no property} \ P \}.$$

The sets have the following meaning

$\mathcal T_P$-the duration of the property $P$;

$\mathcal T_P^0$-the lost duration of the property $P$;

The partition $\{\mathcal T_P, \mathcal T^0_P\}$ of the set
$\mathcal T$ is called $P$ property diagram.

For example, if $P=$commutativity then since any evolution algebra
is commutative, we conclude that any CEA has not commutativity
property transition.

Since an evolution algebra is not a baric algebra, in general, using Theorem
\ref{t2} we can give baric property diagram. Let us do this for the
above given chains $E_i^{[s,t]}$, $i=0,\dots, 24$.

Denote by $\mathcal T^{(i)}_b$ the baric property duration of the CEA $E_i^{[s,t]}$, $i=0,\dots, 24$.
\begin{thm}\label{Te}
 \begin{itemize}
\item[(i)] {\rm (There is no non-baric property transition)} The algebras  $E_i^{[s,t]}$, $i=0,1,2,3,6,10,11,14,22$ are not baric for any time $(s,t)\in \mathcal T$;
\item[(ii)] {\rm (There is no baric property transition)} The algebras  $E_i^{[s,t]}$, $i=16,17,18$ and $E_{23}^{[s,t]}(0,\mu)$, $E_{23}^{[s,t]}(2\mu,\mu)$, $\mu\ne 0$ are baric for any time $(s,t)\in \mathcal T$;
    \item[(iii)] {\rm (There is baric property transition)} The CEAs  $E_i^{[s,t]}$, $i=4,5,7,8,9,12,13$, \, $15,19,20,21,24$ and $E_{23}^{[s,t]}(\lambda,\mu)$, with $\lambda\notin\{0,\mu,2\mu\}$ have baric property transition with baric property duration sets as the following
      $$\mathcal T^{(4)}_b=\left\{(s,t)\in \mathcal T: {\Phi(s)\over\Psi(s)}={\Phi(t)\over\Psi(t)}\right\};\ \
      \mathcal T^{(5)}_b=\left\{(s,t)\in \mathcal T: s\leq t<b, \, \Phi(s)=\Phi(t)\right\};$$
        $$\mathcal T^{(7)}_b=\left\{(s,t)\in \mathcal T: s\leq t<a, \, \Psi(s)=\Psi(t)\right\};\ \
      \mathcal T^{(8)}_b=\left\{(s,t)\in \mathcal T: s\leq t<\min\{a,b\}\right\};$$
        $$\mathcal T^{(9)}_b=\left\{(s,t)\in \mathcal T: t=s+\pi k, k\in \mathbb Z\right\};\ \
      \mathcal T^{(12)}_b=\left\{(s,t)\in \mathcal T: g(s)=\pm{1\over h(s)}\right\};$$ $$
        \mathcal T^{(13)}_b=\left\{(s,t)\in \mathcal T: s\leq t<a, \, \psi(s)=\pm 1\right\};$$ $$
      \mathcal T^{(15)}_b=\left\{(s,t)\in \mathcal T: s\leq t<a, \, \psi(s)=0\right\};\ \
        \mathcal T^{(19)}_b=\left\{(s,t)\in \mathcal T: s\leq t<a\right\};$$ $$
      \mathcal T^{(20)}_b=\left\{(s,t)\in \mathcal T: s\leq t<b\right\}\cup \left\{(s,t)\in \mathcal T: t\geq b, w(s)=0\right\};$$
      $$\mathcal T^{(21)}_b=
      \left\{(s,t)\in \mathcal T: s\leq t<\max\{a,b\}\right\};$$ $$ \mathcal T^{(23)}_b(\lambda,\mu)=
      \left\{(s,t)\in \mathcal T: \theta(t)=\theta(s)\right\}, \lambda\ne 0,\mu,2\mu.$$
      $$\mathcal T^{(24)}_b=
      \left\{(s,t)\in \mathcal T: s\leq t<a\right\}.$$
    \end{itemize}
\end{thm}
\proof By Theorem \ref{t2} a two-dimensional evolution algebra $E^{[s,t]}$ is baric if and only if
  $a_{11}^{[s,t]}\ne 0$, $a_{21}^{[s,t]}=0$ or $a_{22}^{[s,t]}\ne 0$, $a_{12}^{[s,t]}=0$. The assertions of Theorem are
  results of the detailed checking of these conditions.
\endproof

 Note that sets $\mathcal T^{(i)}_b$, $i=8,9,19,21,24$ do not depend on any parameter function. But $\mathcal T^{(i)}_b$, $i=4,5,7,12,13,15,20,23$ depend on some parameter functions and can be controlled by choosing the corresponding parameter functions $\Phi,\Psi,g,h,\psi,w$.  These functions are called {\it baric property
controllers} of the CEAs. Because, they really control the baric
duration set, for example, if some of them is a strong monotone function
then the duration is ``minimal'', i.e. the line $s=t$, but if
a function is a constant function then the baric duration set is
``maximal'', i.e. it is $\mathcal T$. Since these functions are
arbitrary functions, we have a rich class of controller functions,
therefore we have a ``powerful''  control on the property to be
baric.

Now we shall compute the Lebesgue measure $\nu$ of the sets $\mathcal
T^{(i)}_b$. It is easy to see that

$$\nu\left(\mathcal T^{(8)}_b\right)={1\over 2}(\min\{a,b\})^2;\ \ \nu\left(\mathcal T^{(9)}_b\right)=0; \ \ \nu\left(\mathcal T^{(19)}_b\right)={1\over 2}a^2;$$ $$
\nu\left(\mathcal T^{(20)}_b\right)\geq {1\over 2}b^2;\ \
\nu\left(\mathcal T^{(21)}_b\right)={1\over 2}(\max\{a,b\})^2; \ \ \nu\left(\mathcal T^{(24)}_b\right)={1\over 2}a^2.$$

The Lebesgue measure of sets $\mathcal
T^{(i)}_b$, $i=4,5,7,12,13,15,20,23$ depend on the corresponding controller functions.

For a given functions $g$, $h$, $\psi$ one can easily compute $\nu\left(\mathcal T^{(i)}_b\right)$, $i=12,13,15$. For example, if $g$, $h$, $\psi$ are elementary functions which do not have "constant parts" in their graphs then $\nu\left(\mathcal T^{(i)}_b\right)=0$, $i=12,13,15$.

\begin{defn}\label{d6}\cite{clr} A function $\theta$ defined on $\R$ is called a
function of {\it countable variation} if it has the following
properties:

1. it is continuous except at most on a countable set, it has only jump-type
discontinuities;

2.  it has at most a countable set of singular (extremum) points.
\end{defn}

The following theorem gives a characteristics of the baric property
duration set.

\begin{thm}\label{Te4} If the controller function of $\mathcal T^{(i)}_b$,
$i=4,5,7,23$ (i.e. ${\Phi(t)\over \Psi(t)}$ for $i=4$; $\Phi(t)$ for $i=5$; $\Psi(t)$ for $i=7$; $\theta(t)$ for $i=23$) is a function of countable variation,  then the baric duration set $\mathcal T^{(i)}_b$ has zero Lebesgue measure, that is the corresponding CEA is not baric almost surely.
\end{thm}
\proof Is similar to the proof of Theorem 4.8. of \cite{clr}.
\endproof

\section{Absolute nilpotent elements transition.}

Recall that the element $x$ of an algebra $A$ is called an {\it absolute
nilpotent} if $x^2=0$.

Let $E=\R^n$ be an evolution algebra over the field $\R$ with
structural constant coefficients matrix $\M=(a_{ij})$, then for
arbitrary $x=\sum_ix_ie_i$ and $y=\sum_iy_ie_i\in \R^n$ we have
$$xy=\sum_j\left(\sum_ia_{ij}x_iy_i\right)e_j, \ \
x^2=\sum_j\left(\sum_ia_{ij}x^2_i\right)e_j.$$

For a $n$-dimensional evolution algebra $\R^n$ consider operator
$V \colon \R^n\to \R^n$, $x\mapsto V(x)=x'$ defined as
\begin{equation}\label{v}
x'_j= \sum_{i=1}^na_{ij}x_i^2, \ \ j=1,\dots,n.
\end{equation}
This operator is called {\it evolution operator} \cite{ly}.

We have $V(x)=x^2$, hence the equation $V(x)=x^2=0$ is given by the
following system
\begin{equation}\label{n1}
\sum_ia_{ij}x_i^2=0, \ \ j=1,\dots,n.
\end{equation}
 In this section we shall solve the system (\ref{n1}) for $E^{[s,t]}_i$, $i=0,\dots,24.$

For a CEA $E_i^{[s,t]}$ with matrix $\M_i^{[s,t]}$ denote
$$\mathcal T^{(i)}_{nil}=\{(s,t)\in \mathcal T: E_i^{[s,t]}\ \ \mbox{has unique absolute
nilpotent}\}, \ \ \mathcal T_{nil}^0=\mathcal T\setminus \mathcal
T_{nil}.$$

The following theorem gives an answer on problem of existence of
``uniqueness of absolute nilpotent element'' property transition.

\begin{thm}\label{t5} \begin{itemize}

\item[(1)] The CEAs $E_i^{[s,t]}$, $i=3,4,5,9,10,17,22,23,24$ have unique
absolute nilpotent element $(0,0)$ for any time $(s,t)\in \mathcal T$.

\item[(2)] There CEAs $E_i^{[s,t]}$, $i=0,1,2,16,19$  have infinitely many of absolute nilpotent
elements for any time $(s,t)\in \mathcal T$.

\item[(3)] The CEAs $E_i^{[s,t]}$, $i=6,7,8,11,12,13,14,15,18,20,21$ have "uniqueness of
absolute nilpotent element" property transition with the property duration sets as the following
$$\mathcal T^{(i)}_{nil}=\{(s,t)\in \mathcal T: t<a\}, \ \ a>0, \ \ i=6,7,8,11,18.$$
$$\mathcal T^{(12)}_{nil}=\left\{(s,t)\in \mathcal T: g^2(t)\leq {1\over h^2(s)}\right\};$$
$$\mathcal T^{(13)}_{nil}=\left\{(s,t)\in \mathcal T: s\leq t<a, \, \psi^2(s)\leq 1\right\}; \ \
\mathcal T^{(14)}_{nil}=\left\{(s,t)\in \mathcal T: \Phi(s)\psi(s)>0\right\};$$
$$\mathcal T^{(15)}_{nil}=\left\{(s,t)\in \mathcal T: s\leq t<a, \psi(s)>0\right\};$$
$$\mathcal T^{(20)}_{nil}=\left\{(s,t)\in \mathcal T: s\leq t<b\right\}\cup \left\{(s,t)\in \mathcal T: t\geq b,\, {w(s)\over \Phi(s)}>0\right\};$$
$$\mathcal T^{(21)}_{nil}=\left\{(s,t)\in \mathcal T: s\leq t<\min\{a,b\}\right\}\cup \left\{(s,t)\in \mathcal T: b\leq t<a,\, b<a, \, v(s)>0\right\}.$$
\end{itemize}
\end{thm}
\proof The proof consists the simple analysis of the solutions of the
system (\ref{n1}) for each $E^{[s,t]}_i$, $i=0,\dots,24.$
\endproof

\section{Idempotent elements transition}

A element $x$ of an algebra $\A$ is called {\it idempotent} if
$x^2=x$; such points of an evolution algebra are especially
important, because they are the fixed points (i.e. $V(x)=x$) of the
evolution operator $V$, (\ref{v}).  We denote by ${\mathcal Id}(E)$
the set of idempotent elements of an algebra $E$. Using (\ref{v}) the
equation $x^2=x$ can be written as
\begin{equation}\label{v1}
x_j= \sum_{i=1}^na_{ij}x_i^2, \ \ j=1,\dots,n.
\end{equation}
The general analysis of the solutions of the system (\ref{v1}) is
very difficult. We shall solve this problem for the CEAs $E_i^{[s,t]}$,
$i=0,\dots,24$.

The following theorem gives the time-dynamics of the idempotent elements
for algebras  $E_i^{[s,t]}$,
$i=0,\dots,24$.

\begin{thm}\label{ti}\begin{itemize}
\item[(1)] The algebras $E_i^{[s,t]}$, $i=0,1,2,$ have unique idempotent $(0,0)$ any time $(s,t)\in\mathcal T$.
\item[(2)] The algebras $E_i^{[s,t]}$, $i=3,10,12,14,16,22$ have two idempotents $(0,0)$, $(x_i(s,t), y_i(s,t))$ any time $(s,t)\in\mathcal T$. Moreover explicit formula of each $x_i(s,t)$ and $y_i(s,t)$ can be given.
\item[(3)] We have
$${\mathcal Id}\left(E_4^{[s,t]}\right)=\left\{\begin{array}{llll}
\{0,z_1,z_2,z_3\}, \ \ \mbox{if}\ \ (s,t)\in \left\{(s,t)\in\mathcal T: {\Phi(t)\over\Phi(s)}= {\Psi(t)\over\Psi(s)}\right\}\\[2mm]
\{0,z_3\}, \ \ \mbox{if}\ \ (s,t)\in \left\{{\Phi(t)\over\Phi(s)}\ne{\Psi(t)\over\Psi(s)}, \, D(s,t)<0\right\}\\[2mm]
\{0,z_3, (x_*,y_*)\}, \ \ \mbox{if}\ \ (s,t)\in \left\{{\Phi(t)\over\Phi(s)}\ne{\Psi(t)\over\Psi(s)}, \, D(s,t)=0\right\}\\[2mm]
\{0,z_3, (x_{\pm},y_{\pm})\}, \ \ \mbox{if}\ \ (s,t)\in \left\{{\Phi(t)\over\Phi(s)}\ne{\Psi(t)\over\Psi(s)}, \, D(s,t)>0\right\},\\[2mm]
\end{array}\right.$$
where $0=(0,0), \, z_1=(0,{\Phi(t)\over\Phi(s)}),\, z_2=({\Phi(t)\over\Phi(s)},0),\, z_3=({\Phi(t)\over\Phi(s)},{\Phi(t)\over\Phi(s)})$, $D(s,t)={\Phi(t)\over\Phi(s)}\left(2{\Psi(t)\over\Psi(s)}-{\Phi(t)\over\Phi(s)}\right)$. The explicit formulas of $x_*, y_*$, $x_{\pm}$ and $y_{\pm}$ are given below.
The sets $\left\{{\Phi(t)\over\Phi(s)}={\Psi(t)\over\Psi(s)}\right\}$, $\left\{{\Phi(t)\over\Phi(s)}=2{\Psi(t)\over\Psi(s)}\right\}$ are critical (boundary) sets of the idempotent elements transition.
\item[(4)] We have
$${\mathcal Id}\left(E_5^{[s,t]}\right)=\left\{\begin{array}{llll}
\{0,z_1,z_2,z_3\}, \ \ \mbox{if}\ \ (s,t)\in \left\{(s,t)\in\mathcal T: s\leq t<b,\, \Phi(t)=\Phi(s)\right\}\\[2mm]
\{0,z_3\}, \ \ \mbox{if}\ \ (s,t)\in \left\{s\leq t<b, \Phi(t)\ne\Phi(s), \, D(s,t)<0\right\}\\[2mm]
\{0,z_3, (x_*,y_*)\}, \ \ \mbox{if}\ \ (s,t)\in \left\{s\leq t<b, \Phi(t)\ne\Phi(s), \, D(s,t)=0\right\}\\[2mm]
\{0,z_3, (x_{\pm},y_{\pm})\}, \ \ \mbox{if}\ \ (s,t)\in \left\{s\leq t<b,\, \Phi(t)\ne\Phi(s), \, D(s,t)>0\right\},\\[2mm]
\{0, z_3\}, \ \ \mbox{if} \ \ t\geq b.
\end{array}\right.$$
where $z_i$ as in (3), $D(s,t)={\Phi(t)\over\Phi(s)}\left(2-{\Phi(t)\over\Phi(s)}\right)$.
\item[(5)] Algebras $E_i^{[s,t]}$, $i=6,11,13,15,19$ have two idempotent elements for any time $(s,t)$ with $s\leq t<a$ and a unique idempotent for time $(s,t)$ with $t\geq a$. The critical line of the transition is $t=a$.
 \item[(6)] We have
$${\mathcal Id}\left(E_7^{[s,t]}\right)=\left\{\begin{array}{llll}
\{0,z_1,z_2,z_3\}, \ \ \mbox{if}\ \ (s,t)\in \left\{(s,t)\in\mathcal T: s\leq t<a,\, \Psi(t)=\Psi(s)\right\}\\[2mm]
\{0,z_3\}, \ \ \mbox{if}\ \ (s,t)\in \left\{s\leq t<a, \Psi(t)\ne\Psi(s), \, d(s,t)<0\right\},\\[2mm]
\{0,z_3,(x_*,y_*)\}, \ \ \mbox{if}\ \ (s,t)\in \left\{s\leq t<a, \Psi(t)\ne\Psi(s), \, d(s,t)=0\right\},\\[2mm]
\{0,z_3, (x_{\pm},y_{\pm})\}, \ \ \mbox{if}\ \ (s,t)\in \left\{s\leq t<a, \Psi(t)\ne\Psi(s), \, d(s,t)>0\right\},\\[2mm]
0, \ \ \mbox{if} \ \ t\geq a,
\end{array}\right.$$
where $d(s,t)={2\Psi(t)\over\Psi(s)}-1$.
The critical sets are $t=a$, $\Psi(t)=\Psi(s)$, $\Psi(s)=2\Psi(t)$.
 \item[(7)] For $a\leq b$ we have
$${\mathcal Id}\left(E_8^{[s,t]}\right)=\left\{\begin{array}{ll}
\{(0,0),(0,1),(1,0),(1,1)\}, \ \ \mbox{if}\ \ (s,t)\in \left\{(s,t)\in\mathcal T: s\leq t<a\right\}\\[2mm]
(0,0), \ \ \mbox{if} \ \ t\geq a;
\end{array}\right.$$

For $a>b$ we have
$${\mathcal Id}\left(E_8^{[s,t]}\right)=\left\{\begin{array}{lll}
\{(0,0),(0,1),(1,0),(1,1)\}, \ \ \mbox{if}\ \ (s,t)\in \left\{(s,t)\in\mathcal T: s\leq t<b\right\}\\[2mm]
\{(0,0),(1,1)\}, \ \ \mbox{if}\ \ (s,t)\in \left\{(s,t)\in\mathcal T: b\leq t<a\right\}\\[2mm]
(0,0), \ \ \mbox{if} \ \ t\geq a;
\end{array}\right.$$
The lines $t=a$ and $t=b$ are critical for the transition.
\item[(8)] The algebra $E_9^{[s,t]}$, has three idempotent elements $(0,0), (1,0), (0,1)$ for any time $(s,t)$ with $t=s+2\pi n$; has three idempotent elements $(0,0), (-1,0), (0,-1)$ for any time $(s,t)$ with $t=s+(2n+1)\pi$ and at least one  idempotent for time $(s,t)$ with $t\ne s+\pi n$, $n\in \mathbb Z$.
\item[(9)] We have
$${\mathcal Id}\left(E_{17}^{[s,t]}\right)=\left\{\begin{array}{lll}
\{(0,0),z_2\}, \ \ \mbox{if}\ \  D(s,t)<0\\[2mm]
\{(0,0),z_2, ({\Phi(s)\over 2\Phi(t)},{\Psi(s)\over \Psi(t)})\}, \ \ \mbox{if}\ \  D(s,t)=0\\[2mm]
\{(0,0),z_2, (x_{\pm},y_{\pm})\}, \ \ \mbox{if} \ \ D(s,t)>0,
\end{array}\right.$$
where $D(s,t)={4\Phi^2(t)\Psi(s)\over \Phi(s)\Psi^2(t)}(g(t)-g(s))-1$.
\item[(10)] We have
$${\mathcal Id}\left(E_{18}^{[s,t]}\right)=\left\{\begin{array}{llll}
\{(0,0),(1,0)\}, \ \ \mbox{if}\ \  s\leq t<a, \, D(s,t)<0, \\[2mm]
\{(0,0),(1,0), ({1\over 2},{\Psi(s)\over \Psi(t)})\}, \ \ \mbox{if}\ \ s\leq t<a,\, D(s,t)=0,\\[2mm]
\{(0,0),(1,0), (x_{\pm},{\Psi(s)\over\Psi(t)})\}, \ \ \mbox{if} \ \ s\leq t<a,\, D(s,t)>0,\\[2mm]
\{(0,0),({h(t)\Psi(s)\over \Psi^2(t)},{\Psi(s)\over\Psi(t)})\}, \ \ \mbox{if} \ \  t\geq a,\\[2mm]
\end{array}\right.$$
where $D(s,t)=1-{4\Psi(s)(h(t)-h(s))\over\Psi^2(t)}$.
\item[(11)] We have
$${\mathcal Id}\left(E_{20}^{[s,t]}\right)=\left\{\begin{array}{llll}
\{(0,0),z_2\}, \ \ \mbox{if}\ \  s\leq t<b, \, D(s,t)<0, \\[2mm]
\{(0,0),z_2, ({\Phi(s)\over 2\Phi(t)},1)\}, \ \ \mbox{if}\ \ s\leq t<b,\, D(s,t)=0,\\[2mm]
\{(0,0),z_2, (x_{\pm},1)\}, \ \ \mbox{if} \ \ s\leq t<a,\, D(s,t)>0,\\[2mm]
\{(0,0),z_2\} \ \ \mbox{if}, \ \  t\geq b,\\[2mm]
\end{array}\right.$$
where $D(s,t)=1-{4\Phi^2(t)(v(t)-v(s))\over\Phi(s)}$.
\item[(12)] We have
$${\mathcal Id}\left(E_{21}^{[s,t]}\right)=\left\{\begin{array}{llllll}
\{(0,0),(1,0)\}, \ \ \mbox{if}\ \  s\leq t<\min\{a,b\}, \, D(s,t)<0, \\[2mm]
\{(0,0),(1,0),({1\over 2},1)\},  \ \ \mbox{if}\ \ s\leq t<\min\{a,b\}, \, D(s,t)=0,\\[2mm]
\{(0,0),(1,0),(x_{\pm},1)\}, \ \ \mbox{if} \ \ s\leq t<\min\{a,b\}, \, D(s,t)>0,\\[2mm]
\{(0,0),(1,0)\}, \ \ \mbox{if} \ \  b<a, \, b\leq t<a,\\[2mm]
\{(0,0),(v(t),1)\}, \ \ \mbox{if} \ \  b>a, \, a\leq t<b,\\[2mm]
(0,0), \ \ \mbox{if} \ \  t\geq \max\{a,b\},\\[2mm]
\end{array}\right.$$
where $D(s,t)=1-4(v(t)-v(s))$.
\item[(13)] We have
$${\mathcal Id}\left(E_{23}^{[s,t]}(0,\mu)\right)=\left\{\begin{array}{lll}
\{(0,0),(0,1)\}, \ \ \mbox{if}\ \  D(s,t)<0\\[2mm]
\{(0,0),(0,1), ({\theta(s)\over \theta(t)},{1\over 2})\}, \ \ \mbox{if}\ \  D(s,t)=0\\[2mm]
\{(0,0),(0,1), ({\theta(s)\over \theta(t)}, y_{\pm})\}, \ \ \mbox{if} \ \ D(s,t)>0;
\end{array}\right.$$
$${\mathcal Id}\left(E_{23}^{[s,t]}(2\mu,\mu)\right)=\left\{\begin{array}{lll}
\{(0,0),(1,0)\}, \ \ \mbox{if}\ \  D(s,t)<0\\[2mm]
\{(0,0),(1,0), ({1\over 2},{\theta(s)\over \theta(t)})\}, \ \ \mbox{if}\ \  D(s,t)=0\\[2mm]
\{(0,0),(1,0), (x_{\pm}, {\theta(s)\over \theta(t)})\}, \ \ \mbox{if} \ \ D(s,t)>0,
\end{array}\right.$$
where $D(s,t)=1-{4\theta(s)\over \theta(t)}\left({\theta(s)\over \theta(t)}-1\right)$.
 \item[(14)] We have
$${\mathcal Id}\left(E_{24}^{[s,t]}\right)=\left\{\begin{array}{ll}
\{(0,0),(0,1),(1,0),(1,1)\}, \ \ \mbox{if}\ \ (s,t)\in \left\{(s,t)\in\mathcal T: s\leq t<a\right\}\\[2mm]
\{(0,0),\left({g(t)\over (1-g(t))^2+g^2(t)}, {1-g(t)\over (1-g(t))^2+g^2(t)}\right)\},\ \ \mbox{if} \ \ t\geq a.
\end{array}\right.$$
\end{itemize}
\end{thm}
\proof The proof contains detailed analysis of solutions of the system (\ref{v1}) for each $E^{[s,t]}_i$. We shall give here proof of the assertion (3) which is more substantial.
In case of $E^{[s,t]}_4$ the
system (\ref{v1}) has the following form
\begin{equation}\label{v2}\begin{cases}
2x=\left({\Phi(t)\over\Phi(s)}+{\Psi(t)\over\Psi(s)}\right)x^2+\left({\Phi(t)\over\Phi(s)}-{\Psi(t)\over\Psi(s)}\right)y^2;\\
2y=\left({\Phi(t)\over\Phi(s)}-{\Psi(t)\over\Psi(s)}\right)x^2+\left({\Phi(t)\over\Phi(s)}+{\Psi(t)\over\Psi(s)}\right)y^2.
\end{cases}
\end{equation}

 {\it Case} $\Phi(t)\equiv\Psi(t)$. It is easy to see that in this case the
system (\ref{v2}) has only four solutions $0=(0,0)$ and
$$z_1=z_1(s,t)=\left(0,{\Phi(s)\over\Phi(t)}\right), \, z_2=z_2(s,t)=\left({\Phi(s)\over\Phi(t)},0\right),\,
z_3=z_3(s,t)=\left({\Phi(s)\over\Phi(t)}, {\Phi(s)\over\Phi(t)}\right).$$

{\it Case} $\Phi(t)\ne\Psi(t)$. In this case the solutions $0$
and $z_3$ still exist.

{\sl Subcase} $x=0$.  For $x=0$ we have only solution $(0,0)$ if ${\Phi(t)\over\Phi(s)}\ne {\Psi(t)\over\Psi(s)}$ and
there are two solutions $(0,0)$ and $z_1$ if ${\Phi(t)\over\Phi(s)}={\Psi(t)\over\Psi(s)}$.

{\sl Subcase} $y=0$.  For $y=0$ we have only solution $(0,0)$ if ${\Phi(t)\over\Phi(s)}\ne {\Psi(t)\over\Psi(s)}$ and
there are two solutions $(0,0)$ and $z_2$ if ${\Phi(t)\over\Phi(s)}={\Psi(t)\over\Psi(s)}$.

{\sl Subcase} $xy\ne0$. Set $u={x\over y}$.  From system (\ref{v2}) we get

\begin{equation}\label{v3}
(u-1)\left(\left({\Phi(t)\over\Phi(s)}-{\Psi(t)\over\Psi(s)}\right)u^2-2{\Psi(t)\over\Psi(s)}u+\left({\Phi(t)\over\Phi(s)}-{\Psi(t)\over\Psi(s)}\right)\right)=0.
\end{equation}

This equation has unique solution $u=1$ if ${\Phi(t)\over\Phi(s)}={\Psi(t)\over\Psi(s)}$ or ${\Phi(t)\over\Phi(s)}\ne {\Psi(t)\over\Psi(s)}$ and $D(s,t)={\Phi(t)\over\Phi(s)}\left(2{\Psi(t)\over\Psi(s)}-{\Phi(t)\over\Phi(s)}\right)<0$.
For ${\Phi(t)\over\Phi(s)}\ne {\Psi(t)\over\Psi(s)}$ it has two solutions $u=1$ and $u=u_*(s,t)$ if $D(s,t)=0$ and three solutions $u=1$, $u=u_{\pm}(s,t)$ if $D>0$.

Now one can describe $x, y$ corresponding to the solutions of
(\ref{v3}). The case $u=1$, i.e. $x=y$ does not give any new
solution. For $u=u_*, u_{\pm}$ we have $x=u_*y$ and $x=u_{\pm}y$, substituting these in
the second equation of (\ref{v2}) after simple calculations we get
the following non-zero solutions to (\ref{v2}):
$$x_*={\Phi(s)\over\Phi(t)}, \ \ y_*={\Psi(s)\over\Psi(t)}-{\Phi(s)\over\Phi(t)};$$
$$x_{\pm}={1\pm{\Psi(s)\over\Psi(t)}\sqrt{D(s,t)}\over {\Phi(t)\over\Phi(s)}\pm\sqrt{D(s,t)}},
\ \  y_{\pm}={{\Phi(t)\over\Phi(s)}-{\Psi(t)\over\Psi(s)}\over {\Psi(t)\over\Psi(s)}\left({\Phi(t)\over\Phi(s)}\pm\sqrt{D(s,t)}\right)}.$$ Note that
$x_{\pm}, \ \ y_{\pm}$ are well defined for any $(s,t)$ with ${\Phi(t)\over\Phi(s)}\ne {\Psi(t)\over\Psi(s)}$ .
Thus the critical (boundary) times of the transition of idempotent elements are points $(s,t)$ which satisfy ${\Phi(t)\over\Phi(s)}={\Psi(t)\over\Psi(s)}$ or ${\Phi(t)\over\Phi(s)}= {2\Psi(t)\over\Psi(s)}$.

Proofs of the assertions $(i)$, $i=1,2,4-14$ are similar analysis of the system (\ref{v1}) for each $E^{[s,t]}_i$.
\endproof

\section*{ Acknowledgements}

 U. Rozikov thanks Institut des Hautes \'Etudes Scientifiques (IHES), Bures-sur-Yvette, France for
 support of his visit to IHES and he also thanks the Grant No.0251/GF3 of Education and Science Ministry of Republic of Kazakhstan.

{}
\end{document}